\documentclass[a4paper,10pt]{article}

\usepackage{amsmath,amsthm,amssymb}
\usepackage[latin1]{inputenc}
\usepackage[T1]{fontenc}
\usepackage{enumerate}

\theoremstyle{plain}
\newtheorem{thm}{Theorem}[section]
\newtheorem{lem}[thm]{Lemma}

\theoremstyle{definition}

\numberwithin{equation}{section}

\newcommand{\C}{\mathbb C}
\newcommand{\D}{\mathbb D}

\newcommand{\R}{\mathbb R}

\newcommand{\cH}{\mathcal H}
\newcommand{\cO}{\mathcal O}
\newcommand{\id}{\operatorname{id}}

\newcommand{\zbar}{\bar{z}}

\newcommand{\loc}{\mathrm{loc}}

\title{Uniqueness of normalized homeomorphic solutions to nonlinear Beltrami equations}
\date{}

\author{Kari Astala, Albert Clop, Daniel Faraco, \\ Jarmo J\"a\"askel\"ainen, and L\'aszl\'o Sz\'ekelyhidi Jr.}

\begin{document}

\maketitle

\frenchspacing

\begin{abstract}\noindent We settle the problem of the uniqueness  of normalized  homeomorphic solutions to nonlinear Beltrami equations
$\overline\partial f(z)=\cH(z,  \partial f(z))$. It turns out that the uniqueness holds under definite and explicit bounds on the ellipticity at infinity, but not in general.
\end{abstract}

\let\thefootnote\relax\footnotetext{2010 AMS Mathematics Classification Numbers: Primary 35J60; Secondary 30C62.}

\let\thefootnote\relax

\section{Introduction}

Homeomorphic solutions $f \in W^{1,1}_{\loc}(\Omega)$ to the classical Beltrami equation
\begin{equation}\label{bel}
\overline\partial f(z)=\mu(z) \partial f(z), \qquad \| \mu \|_{\infty} \leq k < 1,
\end{equation}
are well-known to be unique up to composing with a conformal mapping. Such solutions
coincide with the class of the two-dimensional quasiconformal mappings, and hence the equation arises naturally in a great variety of topics. For a modern exposition  of the  equation and the quasiconformal mappings in the plane, see the recent monograph  \cite{AIM}. We consider global solutions, solutions in the entire plane $\Omega =\C$. In this case the uniqueness of homeomorphic solutions to \eqref{bel} is obtained simply by requiring that $f(0) = 0$ and $f(1) =1$. We call such homeomorphic solutions $f$ as {\it normalized} solutions to \eqref{bel}.

\medskip

Enquiring the fundamental properties of the {\it nonlinear}  Beltrami equation $\overline\partial f(z)=\cH(z, f(z), \partial f(z))$, the existence of homeomorphic  solutions can be established in great generality. One merely asks of $\cH$ a Lusin type measurability  in the first two variables  and the $k$-Lipschitz condition ($k<1$) in the third; for details, see  Theorem 8.2.1 in \cite{AIM}. The notion of nonlinear Beltrami equations (with more restriction on $\cH$ than above) was introduced in \cite{boj} and \cite{iw}.

However, the uniqueness remains more subtle, even for the system
\begin{equation}
\label{gen}
\overline\partial f(z)=\cH\bigl(z,\partial f(z)\bigr), \qquad \text{for almost every $z\in\C$}.
\end{equation}
In the monograph \cite{AIM} the uniqueness of normalized homeomorphic solutions to  \eqref{gen} was established in the special cases where $\cH(z,w)$ has a compact support in $z$ or when it is homogeneous of degree one in $w$; in particular, we have the uniqueness when  $\cH(z,w)$  is $\R$-linear in $w$. For the general equation \eqref{gen} the question remained open.
\medskip

In this note we show that the uniqueness of normalized homeomorphic solutions holds  if we have small enough bounds on the ellipticity at infinity, but fails in the case of large ellipticity constants. To be more specific,  assume   $\cH:\C\times\C\to\C$  satisfies
\begin{enumerate}
\item[(H1)] For every $w\in\C$, the mapping $z\mapsto \cH(z,w)$ is measurable on $\C$.
\item[(H2)] For $w_1, w_2\in\C$,
$$
|\cH(z,w_1)-\cH(z,w_2)|\leq k(z)|w_1-w_2|, \qquad 0 \leq k(z) \leq k < 1,
$$
for almost every $z\in\C$.
\item[(H3)] $\cH(z,0)\equiv0$.
\end{enumerate}
Our main result is the following.

\begin{thm}\label{uniq}
Suppose $\cH:\C\times\C\to\C$ satisfies (H1)--(H3) for some $k < 1$. If
\begin{equation}\label{atinfty}
\limsup_{|z| \to \infty} k(z) < 3 - 2\sqrt{2} = 0.17157...,
\end{equation}
then the nonlinear Beltrami equation
\begin{equation}\label{gene2}
\overline\partial f(z)=\cH\bigl(z,\partial f(z)\bigr), \qquad \text{for almost every
$z\in\C$},
\end{equation}
admits a unique homeomorphic solution $f \in W^{1,2}_{\loc}( \C )$ normalized by $f(0) = 0$ and $f(1) = 1$.
\medskip

Furthermore, the bound on $k$ is sharp: for each $k > 3 - 2\sqrt{2} $, there are functions
 $\cH:\C\times\C\to\C$ for which  (H1)--(H3) hold, such that \eqref{gene2} admits two  normalized homeomorphic solutions.
\end{thm}

\noindent Note that  in terms of the quasiconformal distortion  the bound \eqref{atinfty} reads as
$$ \limsup_{|z| \to \infty} K(z) < \sqrt{2}, \qquad K(z):=\frac{1+k(z)}{1-k(z)}.
$$

Under extra symmetries in $\cH$ the equation \eqref{gene2}  has a unique normalized solution.
This holds for instance if $\cH(z, tw) \equiv t\cH(z,w)$, no matter
 how large are the ellipticity constants. For another interesting example, note that
the above requirement (H3) asks constant functions to be solutions to the
nonlinear Beltrami equation in question. If we assume, in addition, that also
the identity function satisfies \eqref{gene2} or equivalently
\begin{enumerate}
\item[(H4)] $\cH(z,1) \equiv 0$, 
\end{enumerate}
then ellipticity bounds  slightly weaker
   than  \eqref{atinfty} will suffice:

\begin{thm}\label{uniq1}
Suppose $\cH:\C\times\C\to\C$ satisfies conditions (H1)--(H4) for some $k < 1$. If
\begin{equation}\label{atinfty1}
\limsup_{|z| \to \infty} k(z) < \frac{1}{3},
\end{equation}
then the function $f(z) = z$ is the unique  homeomorphic solution $f \in W^{1,2}_{\loc}( \C )$ to the nonlinear Beltrami equation
\begin{equation}\label{gen3}
\overline\partial f(z)=\cH\bigl(z,\partial f(z)\bigr), \qquad \text{for almost every
$z\in\C$},
\end{equation}
normalized by the conditions  $f(0) = 0$ and $f(1) = 1$.
\medskip

This is complemented with counterexamples:  for any $k > 1/3$ there exists 
 $\cH:\C\times\C\to\C$ satisfying (H1)--(H4)
 such that \eqref{gen3} admits a normalized solution $f \not\equiv z$.
\end{thm}

As it turns out, the knowledge of the existence of enough solutions
gives the uniqueness of normalized solutions. We formulate this as
an abstract theorem and then deduce some corollaries from it.

\begin{thm}\label{flow}
Assume $\cH:\C\times\C\to\C$
satisfies (H1)--(H3) for some $k < 1$. Let $f\in W^{1,2}_{\loc}(\C)$
be a normalized homeomorphic solution to the equation
\begin{equation}\label{gene3}
\overline\partial f(z)=\cH\bigl(z,\partial f(z)\bigr), \qquad
\text{for almost every $z\in\C$}.
\end{equation}
Then $f$ is the unique normalized solution, if  there exists a
continous flow of solutions $\{\psi_t : 0 \leq t  \leq 1\} \subset W^{1,2}_{\loc}(\C)$ of \eqref{gene3}
such that
\begin{itemize}
\item [(F1)] $\psi_0\equiv 0$, $\psi_1=f$,

\item [(F2)] $f-\psi_t$ is quasiconformal, $0 \leq t < 1$,

\item[(F3)]  for fixed $\epsilon > 0$, there exist $R$ and $\delta$ such that 
$\Bigl |\frac{\psi_t(z)-\psi_s(z)}{\psi_t(z)-f(z)}\Bigr| < \epsilon$, when $|z| \geq R$ and $|t - s| < \delta$,   

\item [(F4)] $\psi_t(0)=0$.
\end{itemize}
\end{thm}

Theorem \ref{flow} yields new proofs of the uniqueness
of normalized solutions in some important particular cases; for
instance, when $\cH$ is compactly supported in $z$, the case of the $\R$-linear Beltrami equation or even when $\cH$ is
$1$-homogeneous in $w$, as discussed above. We point out a couple
of further interesting applications. Without the $z$-dependence in
$\cH(z,w)$, every homeomorphic solution is affine.

\begin{thm}\label{nozuniq}
Suppose $\cH:\C\to\C$ is $k$-Lipschitz, $k < 1$, and $\cH(0) = 0$.
Then homeomorphic solutions $f \in W^{1,2}_{\loc}( \C )$ to the nonlinear Beltrami equation
\begin{equation}\label{noz1}
\overline\partial f(z)=\cH\bigl(\partial f(z)\bigr), \qquad \text{for almost every
$z\in\C$},
\end{equation}
are affine; that is, $f(z) = az + \cH(a)\zbar + f(0)$, for some constant $a\in\C$.
\end{thm}

In the case that the identity is a solution, we have the following theorem.

\begin{thm}\label{Lppath}
Suppose that $\cH:\C\times\C\to\C$ satisfies (H1)--(H4) for some $k < 1$. 
If there is a continuous path $\gamma(t):[0,1] \to \C$ such that $\gamma(0)=0, \gamma(1)=1$, and  
uniformly in $t\in[0,1]$
$$ 
\cH\bigl(z,\gamma(t)\bigr) \in L^{p_0}(\C), \qquad \text{for some $p_0 <  2$},
$$
then $f(z)=z$ is the unique normalized $W^{1,2}_{\loc}$-solution  to the nonlinear Beltrami equation
\begin{equation*}
\overline\partial f(z)=\cH\bigl(z,\partial f(z)\bigr), \qquad \text{for almost every
$z\in\C$}.
\end{equation*}
\end{thm}

In particular, if there exists a continuous path of linear solutions connecting $0$ and the identity, then Theorem \ref{Lppath} applies. Nonlinear equations with a rich set of exact solutions enjoy further properties which will be studied in a forthcoming paper.

Finally, we point out an interesting open problem regarding what happens in the borderline case of Theorems \ref{uniq} and \ref{uniq1}. We expect that in this case (i.e., when $\limsup_{|z| \to \infty} k(z) = 3 - 2\sqrt{2}$ or $1/3$, respectively) there is a unique homeomorphic solution $f \in W^{1,2}_{\loc}( \C )$ to the nonlinear Beltrami equation normalized by the conditions  $f(0)=0$ and $f(1) = 1$.

\section{General case, Theorems \ref{uniq} and \ref{uniq1}}

\begin{proof}[Proof of Theorem \ref{uniq}] Let us assume there exist two normalized and homeomorphic solutions $f, g \in W^{1,2}_{\loc}( \C )$ to the nonlinear Beltrami equation \eqref{gene2}. Then conditions (H2) and (H3) imply $|\overline\partial f(z)| \leq k(z) |\partial f(z)|$ and similarly for $g$. Thus $f$, $g$ are quasiconformal. Let
\begin{equation}\label{Kinfty}
K_\infty :=  \limsup_{|z| \to \infty} K(z) < \sqrt{2}, \qquad K(z):=\frac{1+k(z)}{1-k(z)}.
\end{equation}
 Then, for any $K > K_\infty$,
\begin{equation}\label{arvio1}
|f(z)|, |g(z)| \leq C (1 + |z|)^K\!.
\end{equation}
Indeed, we can decompose $f = H \circ F$, where $H$ and $F$ are normalized quasiconformal homeomorphisms with the Beltrami coefficient of $F$  given by  $\chi_{\C \setminus \D(0, R)}\, \mu_{f}$; above $\mu_{f} = \overline\partial f/\partial f$ is the Beltrami coefficient of $f$. Moreover, we may choose  $R$ so large  that $F$ is $K$-quasiconformal in $\C$. Then
$$ \frac{1}{C_K}|z|^{1/K} \leq |F(z)| \leq C_K |z|^K\!, \qquad  |z| \geq 1.
$$
Since $H$ is conformal near $\infty$, $H(z) = cz + \cO(1/z)$, and the bounds \eqref{arvio1} follow.

Next, as $f$, $g$ both satisfy \eqref{gene2}, we have
\begin{equation}\label{erotusqr}
|\overline\partial f(z)  - \overline\partial g(z)| = |\cH\bigl(z, \partial f(z)\bigr) - \cH\bigl(z,  \partial g(z)\bigr)|
\leq k(z)|\partial f(z) - \partial g(z)|,
\end{equation}
for almost every $z \in \C$. Thus the difference is quasiregular, but of course not necessarily injective. By the Sto\"ilow factorization theorem, $f - g = P\circ h$, where $P$ is a holomorphic mapping and $h$ is a normalized $K(z)$-quasiconformal homeomorphism. By \eqref{arvio1} and $\| K \|$-quasiconformality of $h^{-1}$\!, where $\| K \| = \| K\|_{\infty}$, for $|z| \geq 1$,
$$
|P(h(z))| = |f(z) - g(z)| \leq C|z|^K = C|h^{-1}(h(z))|^K  \leq C|h(z)|^{K\| K \|}.
$$
Hence $P$ is a polynomial. Since it has at least two zeroes, points $0$ and $1$, $\deg(P) \geq 2$.

As above, we can decompose $h = H_1 \circ F_1$. Similarly as before: $H_1$ is a normalized quasiconformal mapping and conformal near $\infty$. The mapping $F_1$ is normalized and $K$-quasiconformal in $\C$. This gives us a lower bound for $h$. Combining upper and lower bounds with the fact that $\deg(P) \geq 2$, we achieve, for $|z|$ large enough,
$$
\frac{1}{C}|z|^{2/K} \leq |P(h(z))| = |f(z) - g(z)| \leq C|z|^K\!.
$$
This implies $K \geq \sqrt{2}$ leading to a contradiction with \eqref{Kinfty} when $K > K_\infty$ are sufficiently close.

Our section ''Counterexamples'' below will prove the sharpness of \eqref{atinfty}.
\end{proof}

\bigskip

\begin{proof}[Proof of Theorem \ref{uniq1}] We recall the following topological fact without proof.

\begin{lem}\label{increment}
Let $\gamma$ be a Jordan curve and $f: \C \to \C$ a homeomorphism. Suppose that one of the curves $\gamma$ or $f(\gamma)$ lies inside the other (that is, is separated from $\infty$). Then the increment of the argument
\begin{equation*}
\underset{0\leq t \leq 2\pi}{\Delta}\,\mathrm{arg}\,[f(\xi(t)) - \xi(t)] = \pm 2\pi
\end{equation*}
with the sign depending on the orientation of $f$. Above $\xi$ is any parametrization of $\gamma$.
\end{lem}
\noindent One way to prove the above lemma is to deform the inner curve to a point via a homotopy within the component bounded by the outer curve.

\medskip

Assume now that there exists a normalized solution $\Phi \neq \id$. Conditions (H2) and (H3), and a similar calculation as in \eqref{erotusqr} imply that $\Phi$ and $\Phi-\id$ are $K(z)$-quasiregular, $K(z) = \frac{1+k(z)}{1-k(z)}$.

We have that $\Phi - \id$ is $K(z)$-quasiregular with at least two zeros, points $0$ and $1$. By the Sto\"ilow factorization, $\Phi - \id = P \circ h$, where $P$ is a holomorphic mapping and $h$ is a normalized $K(z)$-quasiconformal homeomorphism. Thus, by the argument principle, for all sufficiently large $R > 0$, the increment of the argument
\begin{equation*}
\underset{|z| = R}{\Delta}\,\mathrm{arg}\,[\Phi(z) - z] \geq 2 \cdot 2\pi.
\end{equation*}
On the other hand, by Lemma \ref{increment}, if the curve $\partial\D(0,R)$ does not intersect the image $\Phi(\partial\D(0,R))$, the increment can be at most $2\pi$. Therefore, for every $R$ large enough, there is a point $z_R$ such that
\begin{equation}\label{zr2}
|\Phi(z_R)| = |z_R| = R.
\end{equation}

The mapping $\Phi$ is a $K(z)$-quasiconformal homeomorphism of the plane and thus \eqref{zr2} forces linear growth at $\infty$. That is, by quasisymmetry,
\begin{equation*}
\frac{1}{\lambda(\|K \|)}|z| \leq |\Phi(z)| \leq \lambda(\|K \|)|z|, \qquad \text{for $|z|$ large enough},
\end{equation*}
where $\|K \| = \| K \|_{\infty}$.
Hence,
\begin{equation}\label{growth}
|\Phi(z) - z| \leq C_{\|K \|}|z|, \qquad \text{for $|z|$ large enough}.
\end{equation}

Similarly as in the proof of Theorem \ref{uniq}, by \eqref{growth} and $\|K \|$-quasiconformality of $h^{-1}$\!,
$P$ is a polynomial. Since it has at least two zeroes, points $0$ and $1$, $\deg(P) \geq 2$.

As before, we can decompose $h = H_1 \circ F_1$, where $H_1$ and $F_1$ are normalized quasiconformal homeomorphisms. Further, $H_1$ is conformal near $\infty$ and $F_1$ is $K$-quasiconformal in $\C$ with $K < 2$. The choice of $K$ can be made by assumption \eqref{atinfty1}. We get a lower bound for $h$. Combining the lower bound with the fact that $\deg(P) \geq 2$ and the upper bound \eqref{growth}, we achieve, for $|z|$ large enough,
$$
c|z|^{2/K} \leq |P(h(z))| = |\Phi(z) - z| \leq C_{\|K \|}|z|.
$$
This is a contradiction, since $K < 2$.

The sharpness is obtained in the next section.
\end{proof}

\section{Counterexamples}

We show that for every $ 3 - 2\sqrt{2} < k < 1$ there is a function  $\cH:\C\times\C\to\C$, measurable in the first variable and satisfying
\begin{equation*}
|\cH(z,w_1)-\cH(z,w_2)|\leq k |w_1-w_2| \qquad \text{and} \qquad \cH(z,0) \equiv 0
\end{equation*}
in the second variable, such that the nonlinear Beltrami equation
\begin{equation}\label{equa}
\overline\partial f(z)=\cH\bigl(z,\partial f(z)\bigr), \qquad \text{for almost every $z\in\C$},
\end{equation}
has at least two different homeomorphic solutions $f \in W^{1,2}_{\loc}(\C)$, normalized by $f(0) = 0$ and $f(1) = 1$.
\medskip

We start the construction by setting, for any $0 < t < 1$,
\begin{align*}
F_t(z) &=
\begin{cases}
(1+t)\,z|z| - t z^2, & \text{for $|z|> 1$,} \\
(1+t)\, z - t z^2, & \text{for $|z| \leq 1$, }
\end{cases} \\
G_t(z) &=
\begin{cases}(1+t)\, z|z| - t z, & \text{for $|z| > 1$,} \\
z, & \text{for $|z| \leq 1$.}
\end{cases}
\end{align*}
Both functions are normalized at $0$ and $1$, and they should be considered as modifications of the radial stretching $\psi(z) = z|z|^{K-1}$, such that their difference   is  a polynomial  vanishing at $0$ and $1$. Hence one may look for a field $\cH(z,w)$ so that $F_t, G_t$ satisfy \eqref{equa}. However, composing with an extra quasiconformal factor we will be able to further reduce the distortion constants. For this purpose take
$$ \varphi(z) = z |z|^{\sqrt{2}-1}, \quad |z| > 1 \quad \mbox{ with }  \quad \varphi(z) = z,  \quad  |z|\leq 1.
$$
and consider the maps $f_t = F_t \circ \varphi^{-1}$ and  $g_t = G_t \circ \varphi^{-1}$, or explicitly 
\begin{align*}
f_t(z) &=
\begin{cases}
(1+t) \,z|z|^{\sqrt{2} -1} - t (z|z|^{1/\sqrt{2} -1})^2, & \text{for $|z|> 1$,} \\
(1+t)\, z - t z^2, & \text{for $|z| \leq 1$, }
\end{cases} \\
g_t(z) &=
\begin{cases}(1+t)\, z|z|^{\sqrt{2} -1} - t z|z|^{1/\sqrt{2} -1}, & \text{for $|z| > 1$,} \\
z, & \text{for $|z| \leq 1$.}
\end{cases}
\end{align*}

Both mappings $f= f_t$ and $g= g_t$ are injective by direct argumentation, and normalized.
It is immediate that
$f - g$ is $K$-quasiregular with $0 < k = \frac{\sqrt{2} - 1}{\sqrt{2} + 1} = 3 - 2\sqrt{2}$, $K = \frac{1+k}{1-k}$.
Directly estimating  $|\overline\partial f(z)|$,  $|\overline\partial g(z)|$ from above and  $|\partial f(z)|$, $|\partial g(z)|$ from below gives that $f$ is $K_f$-quasiregular and $g$ is $K_g$-quasiregular, where
$$
0 < k_f = \frac{\sqrt{2} - 1 + t}{\sqrt{2} + 1 - t} < 1 \quad \text{and} \quad 0 < k_g = \frac{2 - \sqrt{2} + t}{2 + \sqrt{2} + t} < 1.
$$
Next, define for each fixed $z\not\in \partial\D$ the mapping $w \mapsto \cH(z,w)$ as follows. First, fix
\begin{equation}\label{e:cH123}
\cH(z,0)=0,\qquad \cH\bigl(z,\partial f(z)\bigr) = \overline\partial f(z),\qquad \cH\bigl(z,\partial g(z)\bigr)=\overline\partial g(z).
\end{equation}
The computations above show that the map $\cH(z,\cdot):\{0,\partial f(z),\partial g(z)\}\to \C$ is $k_0$-Lipschitz, where $k_0=\max\{k,k_f,k_g\}$. Using the Kirszbraun extension theorem (for example, Theorem 2.10.43 in \cite{Federer}) the mapping can be extended to a $k_0$-Lipschitz map
$\cH(z,\cdot):\C\to \C$. From an abstract use of the Kirszbraun extension theorem, however, it is not entirely clear that the map $\cH$ obtained is measurable in $z$, i.e., that (H1) is satisfied. To show this, one can proceed as follows.

Fix a countable dense set $\mathcal{D}\subset\C$, enumerated as $\mathcal{D}=\{w_4,w_5,w_6,\dots\}$, set $w_1=0$, $w_2=\partial f(z)$, $w_3 = \partial g(z)$, and define $\cH(z,w_k)$ recursively, starting with \eqref{e:cH123}. Assuming $\cH(z,w_k)$ is defined for $k\leq N$ with $N\geq 3$, following \cite{Federer}, we set
$$
Y_s(z)=\bigcap_{j=1}^N\overline{\D}\bigl(a_j(z), s\,r_j\bigr),
$$
where 
$$
a_j(z)=\cH(z,w_j),\quad r_j=k_0|w_j-w_{N+1}|. 
$$
Let
$$
s_0(z):=\inf\{s>0:\,Y_s(z)\neq \emptyset\}.
$$
It is shown in \cite[Lemma 2.10.40]{Federer} that $Y_{s_0}(z)$ consists of a single
point, say $b(z)$, and in the proof of \cite[Theorem 2.10.43]{Federer} that $s_0(z)\le
1$. Furthermore, an additional elementary
argument shows that
$$
(a_1,\dots,a_N)\mapsto b
$$
is a continuous map. Therefore, we set
$$
\cH(z,w_{N+1})=b(z).
$$
Since $a_1(z),a_2(z),a_3(z)$ defined in \eqref{e:cH123} are measurable, it follows recursively that each $a_i(z)$ is measurable in $z$. 
We obtain a $k_0$-Lipschitz map $\cH(z,\cdot):\mathcal{D}\to \C$ such that for each fixed $w\in \mathcal{D}$ the mapping
$z\mapsto \cH(z,w)$ is measurable. Since $\mathcal{D}$ is dense, for each fixed $z$ we can (uniquely) extend $\cH(z,w)$ to a $k_0$-Lipschitz map $\C\to\C$, which
is then measurable in $z$. 

We have now found $\cH(z,w)$, satisfying (H1)--(H3) with $k=k_0$, such that \eqref{equa} has two different normalized solutions. Letting $t\to 0$ makes $k_0 \to 3 - 2\sqrt{2}$. The proof of Theorem \ref{uniq} is thus complete. 
\bigskip

To prove the sharpness  of Theorem \ref{uniq1} one may  modify the counterexample above. However, a more convenient approach is to simply note that given the functions $f_t$ and $ g_t$, one may change the variables so that both the identity and   the composition $\Phi = g_t \circ f_t^{-1}$ satisfy the  same nonlinear Beltrami equation. We may thus use the following  general factorization result to conclude  Theorem \ref{uniq1}.

\begin{lem}\label{Flemma}
Let $\cH:\Omega\times\C\to\C$ be measurable in the first variable, $k(z)$-Lipschitz in the second, $0 \leq k(z) \leq k < 1$, and $\cH(z,0)=0$.
If $f : \Omega \to \Omega'$, $f\in W^{1,2}_{\loc}(\Omega)$, is a homeomorphic solution to the nonlinear Beltrami equation
\begin{equation}\label{Hequation}
\overline\partial f(z)=\cH\bigl(z,\partial f(z)\bigr), \qquad \text{for almost every $z\in\Omega$},
\end{equation}
then any other solution $g \in W^{1,2}_{\loc}(\Omega)$ takes the form $g = \Phi \circ f$, where $\Phi$ solves
\begin{equation}\label{forF}
\overline\partial \Phi(u)=\tilde{\cH}\bigl(u,\partial \Phi(u)\bigr)
\end{equation}
with $\tilde{\cH} : f(\Omega) \times \C \to \C$ measurable in the first variable, $\tilde{k}(u)$-Lipschitz in the second, where $\tilde{k}(u) =\frac{2k(z)}{1+k(z)^2}$, $u =f(z)$, and $\tilde{\cH}(u,0)=0 = \tilde{\cH}(u, 1)$. Furthermore, the function $\tilde{\cH}$ depends only on $\cH$ and the coordinate change $f$, but is independent  of $\Phi$.
\end{lem}

\begin{proof}
We use the chain rule and substitute $g(z) = \Phi(f(z))$ to the equation \eqref{Hequation}. We get the nonlinear relation between $\Phi_u$ and $\Phi_{\bar{u}}$
\begin{equation}\label{nlF}
f_{\zbar}\Phi_u + \overline{f_z}\Phi_{\bar{u}} = \cH\bigl(z, f_z\Phi_u + \overline{f_{\zbar}}\Phi_{\bar{u}}\bigr),
\end{equation}
where $z = f^{-1}(u)$ and $f_{\zbar} = \cH(z, f_z)$. Solving this for $\Phi_{\bar{u}}$ in terms of $\Phi_u$ using the contraction mapping principle, see Chapter 9.1 in \cite{AIM}, gives the equation \eqref{forF}, where $\tilde{\cH} : f(\Omega) \times \C \to \C$ is measurable in the first variable and $\tilde{k}(u)$-Lipschitz in the second.

We are left to check $\tilde{\cH}(u,0)=0 = \tilde{\cH}(u, 1)$. For this we let $\Phi_u = 0$ and $\Phi_u = 1$ in \eqref{nlF} and solve it for $\Phi_{\bar{u}}=\tilde{\cH}(u,0)$ and $\Phi_{\bar{u}}=\tilde{\cH}(u,1)$, respectively. This is equivalent to the equations
\begin{align*}
\overline{f_z}\Phi_{\bar{u}} &= \cH\bigl(z, \overline{f_{\zbar}}\Phi_{\bar{u}}\bigr), \\
f_{\zbar} + \overline{f_z}\Phi_{\bar{u}} &= \cH\bigl(z, f_{z} + \overline{f_{\zbar}}\Phi_{\bar{u}}\bigr).
\end{align*}
In both cases we find that
\[
|\overline{f_z}||\Phi_{\bar{u}}| \leq k^2|\overline{f_z}||\Phi_{\bar{u}}|,
\]
and thus $\Phi_{\bar{u}} = 0$ almost everywhere as wanted. Above we use the $k$-Lipschitz property of $\cH$ and $K$-quasiconformality of $f$, $K = \frac{1+k}{1-k}$, which is a straightforward calculation as in the beginning of the proof of Theorem \ref{uniq}.
\end{proof}

Note that 
$$ k = 3- 2\sqrt{2}\;  \Leftrightarrow \; \frac{2k}{1+k^2} = \frac{1}{3},
$$
thus examples proving sharpness of the bound \eqref{atinfty} yield, via factorization and Lemma \ref{Flemma}, also examples showing the sharpness of Theorem \ref{uniq1}. 
A similar reasoning shows  that the uniqueness part of Theorem \ref{uniq} could be deduced  from Theorem \ref{uniq1}.

\section{Flow of solutions, Theorems \ref{flow}, \ref{nozuniq}, and \ref{Lppath}}

\begin{proof}[Proof of Theorem \ref{flow}]
We use similar methods as in the proof of Theorem 6.2.2 in
\cite{AIM}. Let $f$ be as in the statement of the theorem and $g$ be
another normalized solution.  We construct two different flows of
maps. The flow  $L_t=f-\psi_t$ is a family  of quasiconformal mappings
joining $f$ and $0$. The flow   $g_t=g-\psi_t$ is a family  of
quasiregular mappings joining the homeomorphism $g$ with the noninjective map $g-f$. 

Let $T \subset [0,1)$ denote the set of
parameters $t$ for which $g_{t}$ is a homeomorphism. One such
parameter is $t = 0$. By the Hurwitz-type theorem, Theorem 3.9.4 in
\cite{AIM}, we find that $T$ is a relatively closed subset of
$[0,1)$. Thus we need to show that $T$ is open.

Now, fix a parameter $t\in T$. The mapping
$$g-f = g_t-L_t$$
is, by assumption, a nonconstant $K$-quasiregular mapping with at least two
zeros, points $0$ and $1$. Therefore, the composition
$$(g-f)\circ L_t^{-1}= g_t\circ L_t^{-1}-\id$$
is $K^2$-quasiregular and has also two zeros. 

We use the same ideas as in the proof
of Theorem \ref{uniq1}. First, applying the argument principle and Lemma \ref{increment} to the difference $g_t\circ L_t^{-1}-\id$, we get that for every $R$ large enough there is a point $z_R$ such that
$$|g_t\circ L_t^{-1}(z_R)|=|z_R|=R.$$
As $t\in T, \; g_t\circ L_t^{-1}$ is  quasisymmetric, and since by (F4) it fixes the origin, we obtain  
\begin{equation}\label{growth2}
\frac{1}{\lambda(K)}|z|\leq |g_t\circ L_t^{-1}(z)| \leq \lambda(K)|z| \qquad \text{whenever $|z|\geq R$},
\end{equation}
if $R>0$ is large enough. 

Secondly, the continuity assumption
(F3) and the equation \eqref{growth2} allow us to compare  $g_t \circ
L_t^{-1}$ with $g_s \circ L_t^{-1}$ when $|t-s|$  is small enough.
Indeed, there is $\delta>0$ and $R_0>0$ such that if $|t-s|<\delta$,
one has
$$
|g_s(w)-g_t(w)| = |\psi_t(w) - \psi_s(w)| \leq \frac{1}{2 \lambda(K)}|L_t(w)|,
$$
for $|w|\ge R_0$. Writing $L_t(w)=z$, we obtain
that
\begin{equation}\label{diference}
|g_s\circ L_t^{-1}(z)-g_t\circ L_t^{-1}(z)| \le   \frac{1}{2
\lambda(K)}|z|,
\end{equation}
for every $z$ outside the set $L_t(\D(0,R_0))$ . 

We now fix $w_0\in\C$. Since $g_t\circ L_t^{-1}$ is a
homeomorphism by assumption, the winding number of
$(g_t\circ L_t^{-1})(\partial\D(0,R))$ around $w_0$ is $1$, for $R$ large enough.
Therefore conditions \eqref{growth2} and \eqref{diference} show that for
$|t-s|<\delta$  the winding number of $g_s\circ L_t^{-1}$ is $1$
as soon as $R\ge 2\lambda(K)|w_0|$.  It follows that the
mappings $g_s$  for $|t-s|<\delta$ are homeomorphisms and
$s \in T$. Thus $T$ is open.

We have proven that, for all $t \in [0,1)$, the mappings $g_t$ are
quasiconformal homeomorphisms of the plane. Hence by Hurwitz-type
arguments, e.g.,  Theorem 3.9.4 in \cite{AIM}, their locally uniform limit
$g-f$  is either a homeomorphism or a constant. Having at least two zeroes, it must be the  constant map $0$.
\end{proof}

\medskip

\begin{proof}[Proof of Theorem \ref{nozuniq}]

First we find  a homeomorphic,
 linear, and normalized  solution to the
 nonlinear Beltrami equation \eqref{noz1}: By the Banach fixed point theorem,
 the contraction $w \mapsto 1  - \cH(w)$ has a unique fixed
 point $a\in\C$. Then the  linear mapping $f(z) = a z + \cH(a)\zbar $ is 
 solution to \eqref{noz1}, fixes $0$ and $1$, and is injective by the inequality $|\cH(a)| \leq k|a|$.
 
 Now, we can apply Theorem \ref{flow} with the
 linear maps
 $\psi_t(z) =taz+\cH(ta)\overline{z}$ for $t \in [0,1]$ to see that $f$ is the only normalized solution.
 
 To show that any homeomorphic solution $g$ to \eqref{noz1} is affine, we may assume $g(0)=0$. Given $g(1)=b$, then
 $h(z) = g(z)/b$ is the normalized solution to $\overline\partial h= \widetilde{\cH}\bigl(\partial h)$, where $\widetilde{\cH}(w) = \cH(b \,w)/b$. By the above $h$, and hence $g$, is linear. 
 \end{proof}

\medskip

\begin{proof}[Proof of Theorem \ref{Lppath}]
Since $\cH(z,1)=0$, we already know that $f(z)=z$ is a normalized
solution. We will get the uniqueness by applying Theorem
\ref{flow}. We can assume $\gamma(t) \not\in \{0, 1\}$, when $t \in (0,1)$.

We will construct a concrete flow of solutions. The crucial point is to solve the following nonlinear and inhomogeneous Beltrami equation
\begin{equation}\label{eta}
\overline\partial \eta_t(z)= \cH\bigl(z,\partial\eta_t(z)+\gamma(t)\bigr) \qquad \text{for almost every $z\in\C$}.
\end{equation}
By \cite{AIS}, there exists exactly one solution $\eta_t$ to the above equation \eqref{eta} such that $D \eta_t\in L^p(\C)$.
Namely, this can be established via the invertibility of the nonlinear Beltrami operator $B_t=\mathbf{I}- \cH_t(z, \mathcal{S}) := \mathbf{I} - [\cH(z,\mathcal{S}+\gamma(t))-\cH(z,\gamma(t))]$.  Here  $\mathcal{S}$ stands for the Beurling transform. With this notation, \eqref{eta} gets the form $\overline\partial \eta_t = \cH_t(z,\partial\eta_t) + \cH(z,\gamma(t))$. The $L^p$-invertibility, $1 + k < p < 1 + 1/k$,  of the above operator $B_t=\mathbf{I}- \cH_t(z, \mathcal{S})$   
is proven in \cite{AIS}.  We have assumed that  $\cH(z,\gamma(t))\in L^{p_0}(\C)$ for some $p_0 < 2$. The ellipticity of $\cH$ gives an $L^\infty$-bound,  $|\cH(z,\gamma(t))|\leq k|\gamma(t)|$. Thus have the solution  $\eta_t $ with 
\begin{equation*}
\|\overline \partial \eta_t\|_{L^{{p}}(\C)}\leq C_{{p}}, \qquad \max\{p_0, 1 + k\}  < p < 1 + \frac{1}{k}.
\end{equation*}
 In particular, by the mapping properties of the
Cauchy transform, see, for instance, Theorem 4.3.11 in \cite{AIM}, we obtain a uniform  $L^\infty$-estimate
\begin{equation*}
\| \eta_t\|_{L^{\infty}(\C)}\le C_\infty \qquad \mbox{ with } \;\; \eta_t \in C_0(\hat{\C}).
\end{equation*}

We now claim that $\psi_t(z):=\gamma(t) z+\eta_t(z)-\eta_t(0)$ defines a flow with all the  properties required in  Theorem \ref{flow}. By definition, $\psi_t$ solves the original equation
$$
\overline\partial \psi_t(z)=\cH\bigl(z,\partial \psi_t(z)\bigr) \qquad \text{for almost every $z\in\C$}.
$$
The condition $\cH(z,1)=0$ implies $\eta_1(z) \equiv 0$, and similarly $\eta_0(z) \equiv 0$. Thus $\psi_0(z) = 0$, $\psi_1(z) = z = f(z)$.  Hence we have (F1) and (F4) in Theorem \ref{flow}.  

The quasiconformality of $f-\psi_t$ for $0 \leq t < 1$, condition (F2), follows because
$$
f - \psi_t = \bigl(1 - \gamma(t)\bigr)z - \eta_t(z) + \eta_t(0).
$$
Indeed, since $\eta_t(z) \in C_0(\hat{\C})$, the map  $h:= f-\psi_t$ can be shown  to be
 homotopic to the
homeomorphism $(1-\gamma(t))z + \eta_t(0)$ in $\partial \D(0,R)$ with respect to
$0$, for $R$ large enough. Thus, for example, by Theorem 2.8.1 in \cite{AIM} (for the proof, cf.~Theorem 2.2.4 in \cite{degree}), it follows that $\textrm{deg}(h,0)=1$. Since
$h:\hat{\C} \to \hat{\C}$, the degree is constant and hence equal
to $1$ for each $w \in \hat{\C}$. By quasiregularity of $h$, this
forces $h$ to be a homeomorphism.

It remains to  address the continuity assumption (F3) in Theorem \ref{flow}. Notice that
$$
\overline \partial (\eta_t-\eta_s)=\cH\bigl(z,\partial \eta_t-\gamma(t)\bigr)-\cH\bigl(z,\partial \eta_s-\gamma(s)\bigr).
$$
By applying the Lipschitz condition on $\cH$, we get the pointwise inequality
\begin{equation}\label{estimate}
|\overline \partial (\eta_t-\eta_s)| \le k (|\partial \eta_t-\partial \eta_s|+|\gamma(t)-\gamma(s)|).
\end{equation}
For a compactly supported $\cH$ the continuity estimate follows directly from the invertibility of the Beltrami operators \cite{AIS}, but in  the general case  an additional argument is required. 

Let $\varphi \in C_0^\infty(\D(0,2R))$ with $\varphi(z)=1$, if $|z|\le R$, and $|\nabla \varphi(z)|\le \frac{1}{R}$. Then multiply \eqref{estimate} by $\varphi$  and apply the Caccioppoli type estimate, Theorem 5.4.3 in \cite{AIM}, with the exponent $2<r<1+\frac{1}{k}$ to get
$$\aligned
\biggl(\int_{\D(0,R)}|D(\eta_t-\eta_s)|^r\biggr)^\frac1r
&\leq \biggl( \int_{\D(0,2R)} |\varphi\,D(\eta_t-\eta_s)|^r\biggr)^\frac1r\\
& \le C R^{\frac2r-1}(R |\gamma(t)-\gamma(s)|+ \|\eta_t-\eta_s \|_\infty ).
\endaligned$$ 
Next, we combine this estimate with  the H\"older
estimates for Sobolev functions. More precisely, since by
construction $\eta_t-\eta_s\in W^{1,r}_{\loc}(\C)$ for every
$r\in[2,1+\frac1k)$, at points $z$ with $|z|= R$ one has
\begin{equation}\label{etasetat}
\aligned
 |\eta_t(z)-\eta_s(z)-\eta_t(0)-\eta_s(0)|
 &\leq C\,\biggl(\int_{\D(0,R)}|D(\eta_t-\eta_s)|^r\biggr)^\frac1r\,R^{1-\frac2r}\\
&\leq C\left(|\gamma(t)-\gamma(s)|\,|z|+ C_\infty\right). \endaligned
\end{equation}
Furthermore, the definition of $\psi_t$ and \eqref{etasetat} yield
that there exists a constant $C$ such that if $|z|=R$, then
\begin{equation}\label{psitpsis}
|\psi_t(z)-\psi_s(z)| \le C(|\gamma(t)-\gamma(s)| |z|+ 1).
\end{equation}
On the other hand, clearly $|f(z)-\psi_t(z)|\ge |1-\gamma(t)| |z| - 2 C_\infty$, which for $|z|\ge \frac{ 4 C_\infty}{|1-\gamma(t)|}$ implies that
\begin{equation}
\label{psit}|f(z)-\psi_t(z)|\geq \frac{|1-\gamma(t)|}{2}|z|.
\end{equation}
Combining \eqref{psitpsis} and \eqref{psit} we obtain
\begin{equation*}
\biggl|\frac{\psi_t(z)-\psi_s(z)}{f(z)-\psi_t(z)}\biggr| \le \frac{C}{|1-\gamma(t)|}\biggl(|\gamma(t)-\gamma(s)|+\frac{1}{|z|}\biggr)
\end{equation*}
and the continuity estimate (F3) follows by letting $s\to t$ and $R\to\infty$.
\end{proof}

\textbf{Acknowledgements} Part of the research took place in Madrid where the authors were taken part of the program "Calculus of Variations, Singular Integrals and Incompressible fluids". We would like to thank the warm hospitality of UAM and ICMAT. K.A. was supported by the Academy of Finland, project no. 1134757, the Finnish Centre of Excellence in Analysis and Dynamics Research, and  project MRTN-CT-2006-035651, Acronym CODY, of the European Commission.
A.C. was supported by projects MTM2010-15657 (Spanish Ministry of Science), NF-129254 (Programa Ram\'on y Cajal) and 2009-SGR-420 (Generalitat de Catalunya).  D.F was suported by the spanish grant MTM2008-02568. J.J. was supported by the Academy of Finland, project no. 1134757, and the Vilho, Yrj\"o and Kalle V\"ais\"al\"a Foundation. L.Sz. was supported by the Hausdorff Center for Mathematics in Bonn.

\end{document}